\documentclass[12pt]{amsart}

\usepackage{amsmath}

\usepackage{url}

 \makeatletter
\renewcommand*\subjclass[2][2010]{%
  \def\@subjclass{#2}%
  \@ifundefined{subjclassname@#1}{%
    \ClassWarning{\@classname}{Unknown edition (#1) of Mathematics
      Subject Classification; using '2010'.}%
  }{%
    \@xp\let\@xp\subjclassname\csname subjclassname@#1\endcsname
  }%
}

 \makeatother

\begin{document}

\title[An elemetary proof  of an estimate for a number of primes
\ldots] {An elemetary proof  of an estimate for a number of primes
less  than the product of the first $n$ primes}

\author{Romeo Me\v{s}trovi\'{c}}

{\renewcommand{\thefootnote}{}
  \setcounter{footnote}{0}}
\maketitle

\begin{abstract} 
Let $\alpha$ be a real number such that $1< \alpha <2$ and let 
$x_0=x_0(\alpha)$ be a {\rm(}unique{\rm)} positive solution of the equation 
$$
x^{\alpha-1} -\frac{\pi}{e^2\sqrt{3}}x +1=0.
$$
 Then we prove that for each  positive integer $n>x_0$ there exist at least  
$\left[n^\alpha\right]$ primes between the $(n+1)$th prime and the product of 
the first $n+1$ primes. In particular, we establish a recent Cooke's result 
which asserts that for each positive integer $n$ there are at least 
$n$ primes between the $(n+1)$th prime and the product of the first $n+1$ 
primes. Our proof is based on an  elementary counting method (enumerative 
arguments) and the application of Stirling's formula to give upper 
bound for some binomial coefficients.
   \end{abstract}

Ever since  Euclid of Alexandria, sometimes before
300 {\small B.C.},  first  proved that the number of primes is infinite
(see  Proposition 20 in Book IX  of his legendary {\it Elements}
 \cite{he},
mathematicians have amused themselves by coming up 
with alternate proofs. For more information about the Euclid's proof
of the infinitude of primes see e.g., \cite[p. 414, Ch. XVIII]{d}, 
 and \cite[Section 1]{me}. 
In \cite{me} the author of this article provided a comprehensive 
historical survey of different proofs
of  famous Euclid's theorem on the infinitude of primes
 which has  fascinated generations of mathematicians since its 
first and famous demonstration given by Euclid.
Quite recently, in \cite{me2} the author of this article
presented a  very short and  elementary proof of Euclid's theorem.

Euclid's proof of the infinitude of primes is a paragon of simplicity: 
{\it given a finite list of 
primes $p_1,p_2,\ldots p_n$, multiply them together and add one. 
The resulting number, say $N=p_1p_2\cdots p_n$,
is not divisible by any prime on the list, so any prime factor of $N$ is a
new prime.}

A modification of the above Euclid's proof based on the factorization 
theorem can be found in author's survey article \cite[p. 35, Section 4]{me}. 

Notice that numerous  proofs of the infinitude of primes
yield anyone estimate for distribution of primes \cite{me}.  
Applying  Euclid's proof presented above with $p_1p_2\cdots p_n-1$
instead of $p_1p_2\cdots p_n+1$, we obtain that 
$p_{n+1}<p_1p_2\cdots p_n$ for each $n\ge 2$, where $p_k$ is the $k$th prime.
In 1907 H. Bonse \cite{bo} gave an elementary proof of a stronger inequality, 
now called  Bonse's inequality (for a simple proof   
based on Erd\H{o}s'  method \cite{er} see  \cite[p. 238, Section 4.6]{sa4}):
if $n\ge 4$, then 
  \begin{equation}
p_{n+1}^2<p_1p_2\cdots p_n.
    \end{equation} 
Bonse also proved that $p_{n+1}^3<p_1p_2\cdots p_n$ for all $n\ge 5$.
In 2000 M. Dalezman  \cite[Theorem 1]{dal}
gave an elementary proof of stronger inequality 
$p_{n+1}p_{n+2}<p_1p_2\cdots p_n$ with $n\ge 4$.
 In  1960  L. P\'{o}sa  refined firstly 
 Bonse's inequalities by proving that for every integer $k>1$ 
there is an $n_k$ such that $p_{n+1}^k<p_1p_2\cdots p_n$ for all $n>n_k$.
Further, a syntetic proof (i.e., one not involving the limit concepts of 
analysis) due to S. E.  Mamangakis in 1962 \cite{ma} for a theorem from which 
specialization lead to the following inequalities: $p_{4n}<p_1p_2\cdots p_n$ 
with $n\ge 11$ and $p_{4n}^4<p_1p_2\cdots p_{4n-9}$ with $n\ge 46$. 
In 1971 S. Reich  \cite{re}
showed that for every positive integer $k$ there exists a positive 
integer $n_0=n_0(k)$ such that $p_{n+k}^2<p_1p_2\cdots p_n$ for all
$n\ge n_0$. Furthermore, using a quite different approach from Bonse's, 
in 1988 J. S\'{a}ndor \cite{sa} proved that for $n\ge 3$, 
 $p_1p_2\cdots p_{n-1}+p_n+p_{p_n-2}\le p_1p_2\cdots p_n$
and that for $n\ge 24$,  $p_{n+5}^2+p_{[n/2]}^2<p_1p_2\cdots p_n$   which is 
sharper than Bonse's inequality (1). 

On the other hand, using arguments  of Analytic Number Theory,
many authors have been obtained stronger inequalities than 
those mentioned above. 
In 1977 H. Gupta and S. P. Khare \cite{gk} 
proved that ${n^2\choose n}<p_1p_2\cdots p_n$ for all $n\ge 1794$. 
In 1983 G. Robin \cite[Th\'{e}or\`{e}me 4]{rob} proved that 
$n^n<p_1p_2\cdots p_n$ for each $n\ge 13$. 
Since by Stirling's formula easily follows that 
${n^2\choose n}\sim \frac{e^{n-1}}{\sqrt{2\pi n}}n^n$, it 
follows that the mentioned  inequality by Gupta and  Khare
is stronger than  those due to  Robin.
Motivated by a result of Gupta and Khare, in 2011 H. Alzer and J. S\'{a}ndor 
improved their result \cite[Theorem]{as}.
Moreover, using some Rosser-Schoenfeld's \cite{rs1} and  Robin's estimates 
\cite{rob} for the prime counting function $\pi(x)$
and Chebyshev function $\theta(x)$, in 2000  L. Panaitopol 
\cite{pa1} proved that $p_{n+1}^{n-\pi(n)}<p_1p_2\cdots p_n$ for every 
$n\ge 2$, where $\pi(x)$ is the number of primes $\le x$. 
This improves P\'{o}sa's inequality in the following form: 
$p_{n+1}^k<p_1p_2\cdots p_n$ for $n\ge 2k$ with given $k\ge 1$.
M. Hassani \cite{ha1} refined this inequality in 2006 by proving that for 
$n\ge 101$ the exponent $n-\pi(n)$ can be replaced by 
$(1-1/\log n)(n-\pi(n))$. 
Furthermore, using Panaitopol's  inequality, in  2009
 S.   Zhang showed \cite[Corollary 1]{zh} 
that $2^{p_{n+1}}<p_1p_2\cdots p_n$ for all $n\ge 10$.
This inequality yields an improvement of P\'{o}sa's inequality \cite{pos} 
given above and some Bonse-type inequalities. Furthermore,
various Panaitopol-type  inequalities and related limits are recently 
established by J. S\'{a}ndor \cite{sa2} and  J. S\'{a}ndor 
and A. Verroken \cite{sv}). 
   
Notice that Bonse's inequality, all its refinemenets and improvements 
presented above does not guarantee the existence of  ``many primes" 
less than $p_1p_2\cdots p_n$.  This is also the case with
numerous known elementary proof of Euclid's theorem on the infinitude of 
primes. For example, iterating  the second Bonse's inequality 
we find that $p_{n+2}< (p_1p_2\cdots p_n)^{4/3}$, which repeating
still three  times gives $p_{n+4}\le (p_1p_2\cdots p_n)^{580/729}$. 
This shows that for $n\ge 4$ there exist at least 
$4$ primes between the $n$th prime and the product of the first $n$ primes.
We see from the first Mamangakis'  inequality given above 
that for each  $n\ge 11$ there are at least $3n$ primes  between the $n$th 
prime and $p_1p_2\cdots p_n$. Notice also that the first  S\'{a}ndor's 
inequality presented above and the well known estimation $p_n>n\log n$ with 
$n\ge 5$ (see e.g., \cite[(3.10) in Theorem 3]{rs1}) 
imply that for each $n\ge 2$ there are  at least 
$[n\log n]-2$ primes less than $p_1p_2\cdots p_n$. 
 Quite recently in 2011, applying two simple lemmas in 
the Theory of Finite Abelian Groups related to 
the product of some cyclic groups $\Bbb Z_m$, R. Cooke \cite{c}  
modified Perott's proof from 1881 (\cite{pe}, \cite[page 10]{r2})
to establish that there are at least $n$ primes between the 
$(n+1)$th prime and the product of the first $n+1$ primes. 
Refining the Euler's proof of the infinitude of primes presented below, 
in this note we  improve Cooke's result by proving the following
Bonse-type inequality.

 \vspace{2mm}
\noindent{\bf Theorem 1.} {\it Let $\alpha$ be a real number such that
$1< \alpha <2$ and let $x_0=x_0(\alpha)$ be 
a {\rm(}unique{\rm)} positive solution of the equation 
$$
x^{\alpha-1} -\frac{\pi}{e^2\sqrt{3}}x +1=0.
$$
Then for each  positive integer $n>x_0$ 
there exist at least  $\left[n^\alpha\right]$ primes between 
the $(n+1)$th prime and the product of the first $n+1$ primes.

Moreover,  for each positive integer $n$ there are at least $n$ primes 
between the $(n+1)$th prime and the product of the first $n+1$ primes.}
\vspace{2mm}

\noindent{\bf Remark.} The first assertion of Theorem 1
 can be shortly written in terms of the ``little $o$" notation as
$$
p_{o(n^2)}<p_1p_2\cdots p_{n+1}\leqno(2)
$$ 
as $n\to\infty$, where $p_k$ is the $k$th prime. However, our method applied 
for the proof of Theorem 1 cannot be applied for $\alpha =2$;
namely, this is because of  the inequality (13) with $\alpha =2$ 
is clearly  satisfied for all $m\ge 1$.

A computation via {\tt Mathematica 8} shows that 
$p_{n^2}<(p_n)^2$ for each $n\ge 5$. Namely, it is well known 
(see e.g., \cite{ms})  that $p_n\sim n\log n$ as $n\to\infty$, and so  
$(p_n)^2\sim n^2\log^2 n$ and $p_{n^2}\sim 2n^2\log n$.
This immediately implies that $2(p_n)^2\sim p_{n^2}\log n$.
 Notice also that  using the  known estimates 
$\log n+\log\log n-3/2<p_n/n<\log n+\log\log n-1/2$ with $n\ge 6$
(see e.g., \cite[(3.10) and (3.11) in Theorem 3]{rs1}) and a verification 
via {\tt Mathematica 8} for $1\le n\le 1020$,
easily follows  that $2(p_n)^2> p_{n^2}\log n$ for all $n\ge 1$.
Recall also that it is known that the sequence $(p_n/\log n)$
is strictly increasing (see e.g., \cite[p. 106]{sa3}). 
  \vspace{1mm}

A motivation for our proof given in the next section comes 
from a  less known proof   of Euclid's theorem 
due to Euler  in 1736 (published posthumously in 1862 
\cite{eu1}; also see \cite[Sect. 135]{eu2}, \cite[p. 413]{d}
and cite{me3})  is in fact the first proof of Euclid's theorem after those of Euclid
and C. Goldbach's proof presented in a letter to L. Euler in  July 1730  
(see \cite[p. 6]{r2} and \cite[Appendix C)]{me}).
As noticed in Dickson's History  \cite[p. 413]{d} 
(see also \cite[page 80]{sc}), this proof is also attributed in 
1878/9 by Kummer \cite{k} who gave essentially Euler's argument. 
 This Euler's proof (see e.g., \cite[pp. 134--135]{bu}, \cite{me3}
and \cite[page 3]{pol}; also cf. Pinasco's proof \cite{p}) 
is based on the multiplicativity of the $\varphi$-function 
defined as the number of positive integers not exceeding $n$ and relatively 
prime to $n$. Namely, if $p_1,p_2,\ldots,p_n$ is a list of distinct $n\ge 2$ 
primes with product $P$,  then 
  \begin{equation}
\varphi(P)=(p_1-1)(p_2-1)\cdots (p_n-1)\ge 2^{n-1}\ge 2.
  \end{equation}
The inequality (1) together with the definition of the $\varphi$-function
says there exists at least an integer in the range $[2,P]$ 
that is relatively prime to $P$, but such an integer has a prime 
factor necessarily different from any of the $p_k$ with $k=1,2,\ldots,n$.
This yields the infinitude of primes.

Theorem 1 may be considered as an extension of the above Euler's result 
in order to obtain  the estimate of a number of primes less than 
the product of the first $n$ primes.
Proof of Theorem 1 given in the next section is combinatorial in spirit and 
entirely elementary. It is based on some counting arguments by using 
Stirling's formula to give upper bound for some binomial coefficients.

  \section{Proof of Theorem}

\noindent{\bf Lemma 1.} {\it Let $k$ and $N$ 
be two arbitrary fixed positive integers, and let $N(k,n)$ 
be the number of $k$-tuples  $(x_1,\ldots,x_k)$ of nonnegative 
integers $x_1,\ldots,x_k$ is ${n+k-1\choose k}$ satisfying
the inequality
  $$
\sum_{i=1}^kx_i\le n.
  $$
Then 
 $$
N(k,n)={n+k\choose n}.
 $$}
\begin{proof} 
For a fixed nonnegative integer $m$ with $0\le m\le n$,
denote by $M(k,n)$ the number of $k$-tuples  
$(x_1,\ldots,x_k)$ of nonnegative integers $x_1,\ldots,x_k$ such that
  $$
\sum_{i=1}^kx_i=m.
  $$
Then by induction no $k\ge 1$ it is easy to prove the  well known fact 
that $M(k,m)={m+k-1\choose k-1}$ for each $k\ge 1$. Hence, for such a fixed
$k$ we have
   $$
N(k,n)=\sum_{m=0}^nM(k,m)=\sum_{m=0}^n{m+k-1\choose k-1}.
   $$
Next by induction on $n\ge 0$, using the Pascal's identity
${r\choose i}+{r\choose i+1}={r+1\choose i+1}$ with $0\le i\le r-1$,
we immediately obtain that the sum on the right hand side of the above 
equality is equal to ${n+k\choose k}$. Therefore, we have
$N(k,n)={n+k\choose k}={n+k\choose n}$, as desired.
\end{proof}

\noindent{\bf Lemma 2.} {\it Let $\alpha$ be a real number such that
$\alpha >1$. Then for each positive integer $n$ 
   \begin{equation}
\frac{1}{n!}{[n^\alpha]+n\choose n}<
\frac{e^{2n}(n^{\alpha-1} +1)^n\sqrt{3}}{2 n^{n+1}\pi}
  \end{equation}
where $[a]$ denotes the integer part of $a$. 
Furthermore, for all positive integers $n$ we have
    \begin{equation}
\frac{1}{n!}{2n\choose n}< \frac{2^{2n-1/2}e^n}{ n^{n+1}\pi}.
   \end{equation}
\begin{proof} First observe that (3) holds for $n=1$.
Stirlng's asymptotic formula 
$n!\approx \sqrt{2\pi n}(n/e)^n$ as $n\to\infty$  is 
often presented in the following refined form due to H. Robbins \cite{ro}:
(also see \cite{mi})
   \begin{equation}
n!=\sqrt{2\pi n}\left(\frac{n}{e}\right)^ne^{\gamma_n},
\quad \frac{1}{12n+1}<\gamma_n<\frac{1}{12n},\quad n=1,2,\ldots.
  \end{equation}
Then applying (5)  to all factorials $(\left[n^\alpha\right]+n)!, 
\left[n^\alpha\right]!, n!$, where 
$n^\alpha -1 <a:=\left[n^\alpha\right]\le n^\alpha$
we find that
 \begin{equation}\begin{split}
&\frac{1}{n!}{a+n\choose n}<
\frac{a!}{a!(n!)^2}
 \frac{\sqrt{2\pi (a+n)}\left(
\frac{a+n}{e}\right)^{a+n}}{\sqrt{2\pi a}\left(
\frac{a}{e}\right)^{a}\cdot 
 2\pi n\left(\frac{n}{e}\right)^{2n}}\times\\
& e^{1/(12(a+n))-1/(12a+1)-2/(12n+1)}\\
=&\frac{e^n}{2 n^{2n+1}\pi}\sqrt{1+\frac{n}{a}}\cdot 
\left(1+\frac{n}{a}\right)^{a}(a +n)^n
e^{1/(12(a+n))-1/(12a+1)-2/(12n+1)}.
  \end{split}\end{equation}
Inserting  the inequalities 
 $$
\sqrt{1+\frac{n}{a}}<\sqrt{1+\frac{n}{n^\alpha -1}}
<\sqrt{1+\frac{n}{n-1}}\le\sqrt{3},
  $$
$a+n\le n^\alpha +n$, 
 $$
\left(1+\frac{n}{a}\right)^a= 
\left(\left(1+\frac{n}{a}\right)^{a/n}\right)^n<e^n, 
 $$
and 
$$
\frac{1}{12(a+n)}-\frac{1}{12a+1}-\frac{2}{12n+1}
<\frac{1}{12n}-\frac{2}{12n+1}<0
 $$
into inequality (6), we find that 
   \begin{equation}
\frac{1}{n!}{a+n\choose n}<
\frac{e^{2n}\sqrt{3}(n^\alpha +n)^n}{2\pi n^{2n+1}}=
\frac{e^{2n}(n^{\alpha-1} +1)^n\sqrt{3}}{2 n^{n+1}\pi}.
   \end{equation}
This proves (3). Finally, taking $a=n$ into (6) we immediately 
obtain (4). 
\end{proof}

\noindent{\bf Lemma 3.} {\it The real function $f:[0,+\infty)\to {\mathbf R}$
defined as 
$$
f(x)=e^2\sqrt{3}x^{\alpha-1}-\pi x+1
$$ 
has exactly one positive root. Moreover, if $x_0=x_0(\alpha)$ is this root, 
then $f$
decreases on $[x_0,+\infty)$.} 

\begin{proof}
The  derivative of the function $f$ is 
$f'(x)=e^2\sqrt{3}(\alpha-1)x^{\alpha-2}-\pi$ whose a unique real root
is 
$x_1=\left(e^2\sqrt{3}(\alpha -1)/\pi\right)^{1/(2-\alpha)}$.
Therefore, in view of the fact that  $\alpha-2<0$, we infer that
$f$ is decreasing on  $[x_1,\infty)$. However, since 
$f(x_1)/x_1=\pi(2-\alpha)/(\alpha -1)>0$ and so $f(x_1)>0$, 
and $\lim_{x\to +\infty}f(x)=-\infty$, we conclude that there 
exists a unique positive root $x_0$ of the function $f$.
As $x_0>x_1$ we see that $f$ decreases on $[x_0,+\infty)$.  
 \end{proof}
  \begin{proof}[Proof of Theorem.] First  consider the case 
when $1<\alpha <2$. Then suppose that the assertion is not true.
This means that there are  
an $\alpha$  with $1<\alpha <2$ and 
a  positive integer $m>x_0=x_0(\alpha)$ for which 
there are less  than $\left[m^\alpha\right]$ primes between 
the $(m+1)$th prime and the product of the first $m+1$ primes.
For such a $m$, let $2=p_1<3=p_2<\cdots <p_{m+1}$
be first $m+1$ consecutive primes. Accordingly, suppose that 
$p_{m+2},\ldots,p_{m+1+k}$ are all the primes between $p_{m+1}$ and 
the product $P:=p_1p_2\cdots p_{m+1}$ with $k\le \left[m^\alpha\right]-1$. 
Then  every positive integer less than $P$ and relatively prime to $P$ can be 
factorized  as $p_{m+2}^{x_1}\cdots p_{m+k+1}^{x_k}$ with 
nonnegative integers $x_1,\ldots ,x_k$. Then obviously
we have 
  \begin{equation*}
p_{m+1}^{x_1+\cdots +x_k}<p_{m+2}^{x_1}\cdots p_{m+1+k}^{x_k}<P=
p_1p_2\cdots p_{m+1}< p_{m+1}^{m+1},
  \end{equation*}
whence it follows that 
 \begin{equation}
x_1+\cdots +x_k\le m.
  \end{equation}
Then by Lemma 1, a number $N(k, m)$ of $k$-tuples  
$(x_1,\ldots,x_k)$ of nonnegative integers $x_i$
satisfying the inequality (7) is equal to ${m+k\choose k}$. 
It follows that must be 
   \begin{equation}
{m+k\choose m}\ge 
\varphi(P)=(p_1-1)(p_2-1)\cdots (p_{m+1}-1).
  \end{equation}
On the other hand, since  $p_i-1\ge 2i$ for all $i\ge 5$, it follows that  
for every $n\ge 3$
    \begin{equation}\begin{split}
&(p_1-1)(p_2-1)(p_3-1)(p_4-1)\cdots (p_{n+1}-1)\\
&\ge 48\cdot  2^{n-3}\cdot 5\cdot 6\cdots (n+1)=2^{n-2}(n+1)!.
\end{split}\end{equation}
Then (9) and (10) with $n=m\ge 3$ yield 
 \begin{equation}
 2^{m-2}(m+1)\le \frac{1}{m!}{m+k\choose m}.
  \end{equation}
Using the fact that the sequence $k\mapsto {m+k\choose m}$ 
($k=1,2,\ldots$) is increasing, the previous assumption 
$k<\left[m^\alpha\right]$ and the inequality (3) of Lemma 2, 
we obtain that if  $m\ge 3$, then 
  \begin{equation}
\frac{1}{m!}{m+k\choose m}<\frac{1}{m!}{[m^\alpha]+m\choose m}<
\frac{e^{2m}(m^{\alpha-1} +1)^m\sqrt{3}}{2 m^{m+1}\pi}.
  \end{equation}
Now from (11) and (12) it follows that if $m\ge 3$ then
  \begin{equation*}
2^{m-2}(m+1)<\frac{e^{2m}(m^{\alpha-1} +1)^m\sqrt{3}}{2 m^{m+1}\pi},
  \end{equation*}
or equivalently,
 \begin{equation}
m^{\alpha-1} +1>\frac{2}{e^2}\sqrt[m]{\frac{m(m+1)\pi}{2\sqrt{3}}}\cdot m.
  \end{equation}
Notice that the inequality (13) also holds for $m=1$ and $m=2$.

Since $\pi/(2\sqrt{3})<1$, it follows that  
    \begin{equation*}
\sqrt[m]{\frac{m(m+1)\pi}{2\sqrt{3}}}>
\sqrt[m]{\frac{\pi}{2\sqrt{3}}} 
\ge \frac{\pi}{2\sqrt{3}}.
  \end{equation*}
Substituting this into (13) we have 
 \begin{equation}
m^{\alpha-1} +1-\pi m/(e^2\sqrt{3}) >0.
 \end{equation}
However, using the notations of  Lemma 3, since $m>x_0$ 
this lemma gives
 \begin{equation*}
m^{\alpha-1} +1-\pi m/(e^2\sqrt{3})=f(m)<f(x_0)=0.
 \end{equation*} 
This contradicts (14), and hence the proof when $1<\alpha <2$ 
is finished.

It remains to prove the assertion for $\alpha =1$. 
Then as in the previous case, suppose that for some $m\ge 1$ there exist 
less  than $m$ primes between $p_{m+1}$ and the product 
$p_1p_2\cdots p_{m+1}$. Then as in the previous case with the
condition $k\le m-1$ instead of $k\le \left[m^{\alpha}\right]-1$,
we arrive to the following inequality  analogous to (11):
 \begin{equation}
 2^{m-2}(m+1)\le \frac{1}{m!}{2m-1\choose m}=\frac{1}{2}\cdot 
\frac{1}{m!}{2m\choose m}, \quad m=1,2,\ldots.
  \end{equation}
Next from (15) and the inequality (4) of Lemma 2 for all
$m\ge 1$ we get 
 \begin{equation*}
 2^{m-2}(m+1)<\frac{2^{2m-3/2}e^m}{m^{m+1}\pi},
  \end{equation*}
whence it follows that if $m\ge 6$  then
 \begin{equation*}
1<\frac{2^{m+1/2}e^m}{m^{m+1}(m+1)\pi}=\frac{\sqrt{2}}{m(m+1)\pi}
\left(\frac{2e}{m}\right)^m<1.
  \end{equation*}
A contradiction, and hence our assertion is true if $m\ge 6$.
We immediately verify that  between $p_{n+1}$ and the product 
$p_1p_2\cdots p_{n+1}$ there are at least $n$ primes for all
$n\in\{1,2,3,4,5\}$. This completes the proof.
 \end{proof}

\vspace{5mm}

\noindent{\scriptsize {\rm DEPARTMENT OF MATHEMATICS, MARITIME FACULTY KOTOR,

\noindent UNIVERSITY OF MONTENEGRO

\vspace{-1mm}

\noindent DOBROTA 36, 85330 KOTOR, MONTENEGRO}}

\vspace{-1mm}
\noindent{\it E-mail address}: {\tt romeo@ac.me}

\end{document}